\documentclass[11pt,a4paper]{article}

\usepackage{epsf,epsfig,amsfonts,amsgen,amsmath,amstext,amsbsy,amsopn,amsthm
}
\usepackage{amsmath,times,mathptmx}
\usepackage{amsfonts,amsthm,amssymb}
\usepackage{amsfonts}
\usepackage{graphics}
\usepackage{latexsym,bm}
\usepackage{amsfonts,amsthm,amssymb,bbding}
\usepackage{indentfirst}
\usepackage{graphicx}
\usepackage{color}
\usepackage[colorlinks=true,anchorcolor=blue,filecolor=blue,linkcolor=blue,urlcolor=blue,citecolor=blue]{hyperref}
\usepackage{float}
\usepackage{tikz}
\newtheorem{thm}{Theorem}

\newtheorem{lemma}{Lemma}
\newtheorem{false statement}{False statement}

\theoremstyle{definition}

\newtheorem{claim}{Claim}

\newtheorem{case}{Case}
\newtheorem{subcase}{Subcase}[case]

\newtheorem{problem}{Problem}

\begin{document}
\allowdisplaybreaks[4]
\title{Extremal spectral radius and $g$-good $r$-component connectivity
\footnote{Supported by NSFC (Nos. 12361071 and 11901498).}}
\author{{Wenxiu Ding, Dan Li}\thanks{Corresponding author. E-mail: ldxjedu@163.com.}{, Yu Wang}\\
{\footnotesize College of Mathematics and System Science, Xinjiang University, Urumqi 830046, China}}
\date{}

\maketitle {\flushleft\large\bf Abstract:}
For $F\subseteq V(G)$, if $G-F$ is a disconnected graph with at least $r$ components and each vertex $v\in V(G)\backslash F$ has at least $g$ neighbors, then $F$ is called a $g$-good $r$-component cut of $G$. The $g$-good $r$-component connectivity of $G$, denoted by $c\kappa_{g,r}(G)$, is the minimum cardinality of $g$-good $r$-component cuts of $G$. Let $\mathcal{G}_n^{k,\delta}$ be the set of graphs of order $n$ with minimum degree $\delta$ and $g$-good $r$-component connectivity $c\kappa_{g,r}(G)=k$. In the paper, we determine the extremal graphs attaining the maximum spectral radii among all graphs in $\mathcal{G}_n^{k,\delta}$. A subset $F\subseteq V(G)$ is called a $g$-good neighbor cut of $G$ if $G-F$ is disconnected and each vertex $v\in V(G)\backslash F$ has at least $g$ neighbors. The $g$-good neighbor connectivity $\kappa_g(G)$ of a graph $G$ is the minimum cardinality of $g$-good neighbor cuts of $G$. The condition of $g$-good neighbor connectivity is weaker than that of $g$-good $r$-component connectivity, and there is no requirement on the number of components. As a counterpart, we also study similar problem for $g$-good neighbor connectivity.

\vspace{0.1cm}
\begin{flushleft}
\end{flushleft}

\section{Introduction}
In the paper, we only consider simple connected graphs. Let $G=(V(G), E(G))$ be a graph with vertex set $V(G)$ of order $n=|V(G)|$ and edge set $E(G)$. The adjacency matrix $A(G)$ of $G$ is an $n\times n$ matrix whose $(i,j)$ entry is 1 if $v_{i}$ is adjacent to $v_{j}$ in $G$, and 0 otherwise. It is obvious that $A(G)$ is a real symmetric matrix. Thus its eigenvalues are real numbers. The largest eigenvalue of $A(G)$ is called the spectral radius of $G$, denoted by $\rho(G)$.
Let $N_G(v)$ be the set of neighbors of $v$ in $G$ and $N_G[v]=N_G(v)\cup\{v\}$. The degree of the vertex $v$ is $d_G(v)=|N_G(v)|$. A graph $G$ is called $t$-regular if every vertex has the same degree equal to $t$. Denote by $\delta(G)$ (or for short $\delta$) and $\Delta(G)$ the minimum and maximum degrees of the vertices of $G$. We know $\delta(G)\leq \rho(G)\leq \Delta(G)$, and the equality holds in either of these inequalities if and only if $G$ is regular.

Connectivity is a fundamental concept in graph theory. For a set $F\subseteq V(G)$ (or $F\subseteq E(G)$), the notation $G-F$ represents the graph obtained by removing the vertices (or edges) in $F$ from $G$. If $G-F$
is disconnected, then $F$ is called a vertex (or edge) cut. The classic connectivity (connectivity
for short) of $G$, denoted by $\kappa(G)$, is defined as the cardinality $|F|$ of the smallest set $F\subseteq V(G)$ such that $G-F$ is either disconnected or a graph with a single vertex. The edge connectivity of a graph $G$, denoted by $\lambda(G)$, is the minimal number of edges whose removal produces a disconnected graph.
The study of connectivity has been shown to be very important to graphs and has many applications in measuring reliability and fault-tolerance networks, such as \cite{F. Cao, M.J. Ma}. In general, the larger $\kappa(G)$ or $\lambda(G)$ is, the more reliable the network is. It is well known that $\kappa(G)\leq \lambda(G)\leq \delta(G)$. A graph $G$ is called maximally edge connected or $\lambda$-optimal if $\lambda(G)=\delta(G)$ and maximally vertex connected if $\kappa(G)=\delta(G)$. However, the classic connectivity has certain limitations
and underestimates the resilience of the network. To measure the fault tolerance of an interconnection network
more accurately, Harary \cite{F. Harary} introduced the concept of the conditional connectivity by placing some requirements on the components of $G-F$, where $F$ is a subset of edges or vertices. In 1984, Chartrand, Kapoor, Lesniak and Lick \cite{G. Chartrand} introduced the $r$-component connectivity, as a generalization of the classic connectivity. The $r$-component connectivity $c\kappa_r(G)$ of a non-complete graph $G$ is the minimum number of vertices whose deletion results in a disconnected graph with at least $r$ components or a graph with fewer than $r$ vertices. Clearly, $c\kappa_2(G)=\kappa(G)$. In addition, Latifi, Hedge and Naraghi-Pour \cite{S. Latifi} introduced the concept of $g$-good neighbor connectivity. A subset $F\subseteq V(G)$ is called a $g$-good neighbor cut of $G$ if $G-F$ is disconnected and each vertex $v\in V(G)\backslash F$ has at least $g$ neighbors. The $g$-good neighbor connectivity of $G$, denoted by $\kappa_g(G)$, is the cardinality of a minimum $g$-good neighbor cut of $G$.
However, the above mentioned conditional connectivity is only a condition that restricts the components of $G-F$. In some cases, we need to have various restrictions on the components of $G-F$.
Recently, Zhu, Zhang, Zou and Ye \cite{B. Zhu} gave the concept of the $g$-good $r$-component connectivity. For $F\subseteq V(G)$, if $G-F$ is disconnected and there are at least $r$ components and each vertex $v\in V(G)-F$ has at least $g$ neighbors, then $F$ is called a $g$-good $r$-component cut of $G$. The $g$-good $r$-component connectivity of $G$, denoted by $c\kappa_{g,r}(G)$, is the minimum cardinality of $g$-good $r$-component cuts of $G$. As a significant tool, $g$-good $r$-component connectivity has improved the accuracy of reliability and fault tolerance analysis of networks.

Brualdi and Solheid \cite{R.A. Brualdi} proposed the following general problem, which became one of the classical problems of spectral graph theory.
\begin{problem} \label{problem1} \cite{R.A. Brualdi}
Given a set $\mathcal{G}$ of graphs, find min$\{\rho(G): G\in \mathcal{G}\}$ and max$\{\rho(G): G\in \mathcal{G}\}$, and characterize the corresponding extremal graphs.
\end{problem}

A lot of results concerning the Brualdi-Solheid problem were presented, and some of these results were exhibited in a recent monograph on the spectral radius by Stevanovi$\acute{c}$ \cite{D. Stevanovic}.
This paper focuses on the research connectivity from spectral perspectives. The spectral conditions for the connectivity of graphs have been well investigated. Berman and Zhang \cite{A. Berman} studied the spectral radius of graphs with $n$ vertices and $k$ cut vertices and described the graph that has the maximal spectral radius in this class. The union $G_1\cup G_2$ is defined to be $G_1\cup G_2=(V_1\cup V_2,~E_1\cup E_2)$. The join $G_1\vee G_2$ is obtained from $G_1\cup G_2$ by adding all the edges joining a vertex of $G_1$ to a vertex of $G_2$. Denote by $K_n$ the complete graph with order $n$. $K_k\vee(K_1\cup K_{n-k-1})$ is shown to be the graph with the maximal spectral radius among all graphs of order $n$ with connectivity $\kappa(G)\leq k$ in \cite{J. Li, M.L. Ye}. Lu and Lin \cite{H.L. Lu} proved that $K_k\vee(K_{\delta-k+1}\cup K_{n-\delta-1})$ is the graph with the maximum spectral radius among all graphs of order $n$ with $\kappa(G)\leq k\leq \delta(G)$. Recently, Fan, Gu and Lin \cite{D.D. Fan3} extended some results on classic connectivity. They determined the graphs with maximal spectral radius among all graphs of order $n$ with given minimum degree $\delta$ and $r$-component connectivity. On the other hand,
the spectral conditions for the connectivity of digraphs have received a lot of attention of researchers, see \cite{D.D. Fan3, H.Q. Lin}. Moreover, the spectral conditions for the edge-connectivity of graphs have also been well investigated by many researchers, see \cite{D.D. Fan2, W.J. Ning, Y. Wang}.

Since the $g$-good $r$-component connectivity is consider as a generalization of the connectivity, based on Problem \ref{problem1}, motivated by the above spectral results, we are interested in studying the spectral condition for $g$-good $r$-component connectivity.
\begin{problem}\label{problem2}
Which graphs attain the maximum spectral radii among all connected graphs of order $n$ with fixed minimum degree $\delta$ and $g$-good $r$-component connectivity $c\kappa_{g,r}(G)$?
\end{problem}

For convenience, we use $k$ instead of $c\kappa_{g,r}(G)$. Let $\mathcal{G}_n^{k,\delta}$ be the set of graphs of order $n$ with minimum degree $\delta$ and $g$-good $r$-component connectivity $k$.
Note that there is no bound between $k$ and $\delta$.
Let $G_{n,(g+1)^{r-1}}^{\delta,0}$ be the graph obtained from $K_1\cup (K_{k-1}\vee(K_{n-(r-1)(g+1)-k}\cup (r-1)K_{g+1}))$ by adding $\delta$ edges between the isolated vertex $K_1$ and $K_{k-1}$.
Let $G_{n,(g+1)^{r-1}}^{\delta-g,g}$ be the graph obtained from $K_1\cup (K_k\vee(K_{n-(r-1)(g+1)-k}\cup (r-2)K_{g+1}\cup K_g))$ by adding $\delta-g$ edges between the isolated vertex $K_1$ and $K_k$, and adding $g$ edges between $K_1$ and $K_g$.
Let $G_{n,(g+1)^{r-1}}^{k-1,\delta-k+1}$ be the graph obtained from $K_1\cup (K_{k-1}\vee(K_{n-(r-1)(g+1)-k}\cup (r-1)K_{g+1}))$ by adding $k-1$ edges between the isolated vertex $K_1$ and $K_{k-1}$, and adding $\delta-k+1$ edges between $K_1$ and $K_{n-(r-1)(g+1)-k}$.
Let $G_{n,(g+1)^{r-1}}^{0,\delta}$ be the graph obtained from $K_1\cup (K_{n-(r-1)(g+1)-1}\cup (r-1)K_{g+1})$ by adding $\delta-r+1$ edges between the isolated vertex $K_1$ and $K_{n-(r-1)(g+1)-1}$, and adding one edge between $K_1$ and each $K_{g+1}$, respectively (see Figure \ref{fig1}).
In this paper, we mainly give the answer to Problem \ref{problem2}. Then we have the following result.

\begin{figure}[H] \label{fig1}
\begin{center}
\includegraphics[width=14.5cm, height=8.5cm]{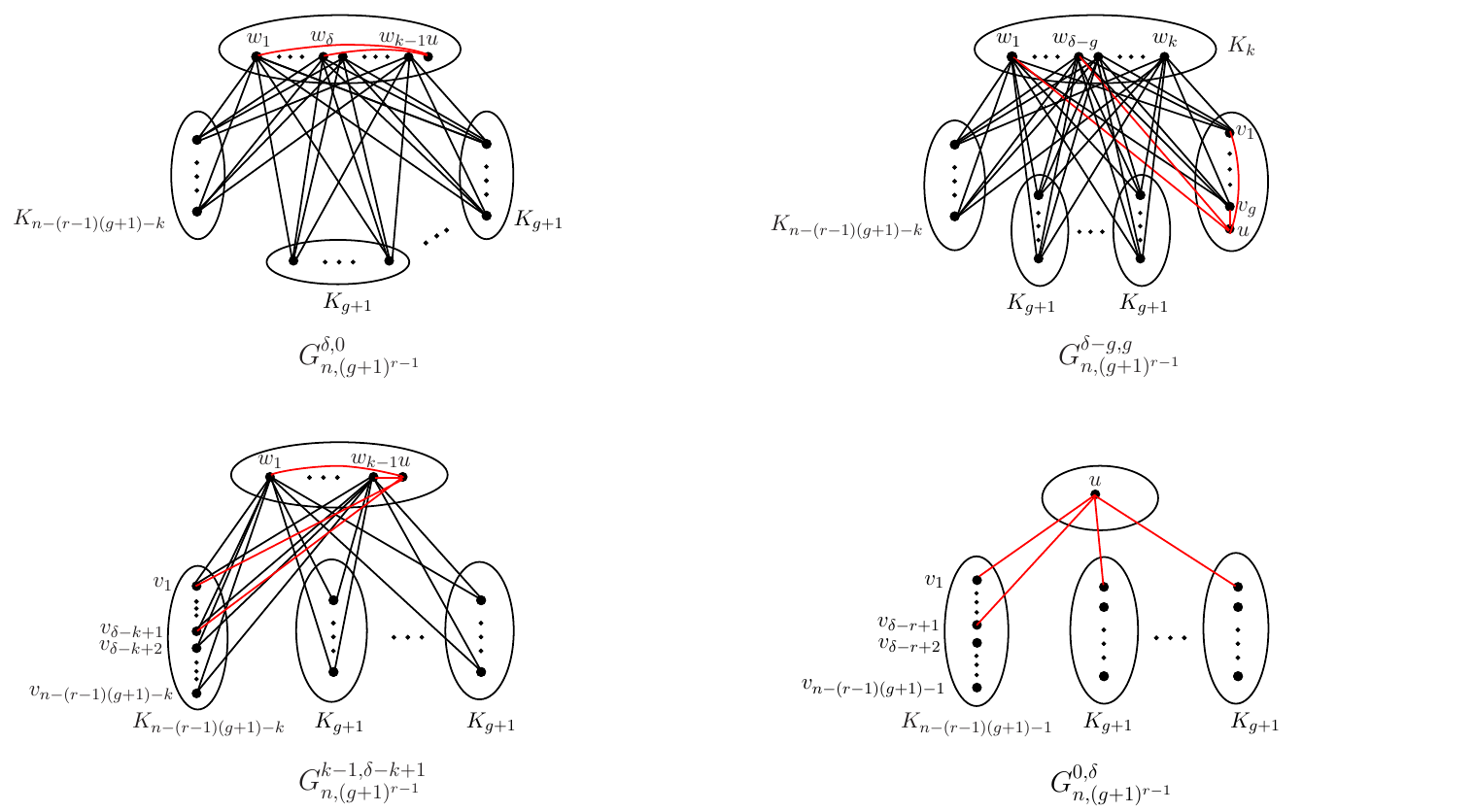}
\caption{$G_{n,(g+1)^{r-1}}^{\delta,0},~G_{n,(g+1)^{r-1}}^{\delta-g,g},~G_{n,(g+1)^{r-1}}^{k-1,\delta-k+1}~\text{and}~ G_{n,(g+1)^{r-1}}^{0,\delta}$}
\label{fig1}
\end{center}
\end{figure}

\begin{thm} \label{thm1}
Let $G\in \mathcal{G}_n^{k,\delta}$, where $n\geq k+r(g+1)$. Then we have the following statements.
\begin{description}
   \item[(I)] If $k>\delta$ and $\delta< g$, then $\rho(G)\leq \rho (G_{n,(g+1)^{r-1}}^{\delta,0})$ with equality if and only if $G\cong G_{n,(g+1)^{r-1}}^{\delta,0}$;
   \item[(II)] If $k>\delta\geq g$, then $\rho(G)\leq \rho(G_{n,(g+1)^{r-1}}^{\delta-g,g})$ with equality if and only if $G\cong G_{n,(g+1)^{r-1}}^{\delta-g,g}$;
   \item[(III)] If $2\leq k\leq\delta<g$, then $\rho(G)\leq \rho(G_{n,(g+1)^{r-1}}^{k-1,\delta-k+1})$ with equality if and only if $G\cong G_{n,(g+1)^{r-1}}^{k-1,\delta-k+1}$;
   \item[(IV)] If $1=k\leq\delta<g$, then $\rho(G)\leq \rho(G_{n,(g+1)^{r-1}}^{0,\delta})$ with equality if and only if $G\cong G_{n,(g+1)^{r-1}}^{0,\delta}$;
   \item[(V)] If $k\leq\delta$ and $g\leq \delta <g+k$, then $\rho(G)\leq \rho(G_{n,(g+1)^{r-1}}^{\delta-g,g})$ with equality if and only if $G\cong G_{n,(g+1)^{r-1}}^{\delta-g,g}$;
   \item[(VI)] If $\delta \geq g+k$, then $\rho(G)\leq \rho(K_k\vee(K_{n-k-(\delta-k+1)(r-1)}\cup (r-1)K_{\delta-k+1}))$ with equality if and only if $G\cong K_k\vee(K_{n-k-(\delta-k+1)(r-1)}\cup (r-1)K_{\delta-k+1})$.
 \end{description}
\end{thm}

The rest of the paper is structured as follows. In Section \ref{sec2}, we determine the graphs attaining the maximum spectral radii among all connected graphs with given minimum degree and $g$-good $r$-component connectivity. In Section \ref{sec3}, we extend the result of extremal spectral radius and characterize the graphs with maximal spectral radii among all connected graphs with given minimum degree and $g$-good neighbor connectivity.

\section{Proof of Theorem \ref{thm1}} \label{sec2}
In this section, we present the proof of Theorem \ref{thm1}. First of all, we list several lemmas that will be used in the follows.
Note that for a connected graph $G$, there is a positive unit eigenvector $\bm{x}$ of $A(G)$ corresponding to $\rho(G)$. Such a unit eigenvector is also called Perron vector of $A(G)$. The component of $\bm{x}$ corresponding to a vertex $v\in V(G)$ is usually written as $x(v)$.

\begin{lemma} \cite{W.J. Ning} \label{lem1}
Let $u,v$ be two distinct vertices of a connected graph $G$, and let $\bm{x}$ be the Perron vector of $A(G)$.
\begin{description}
   \item[(1)] If $ N_G(v)\backslash\{u\}\subset N_G(u)\backslash \{v\}$, then $x(u)>x(v)$;
   \item[(2)] If $N_G(v)\subseteq N_G[u]$ and $N_G(u)\subseteq N_G[v]$, then $x(u)=x(v)$.
 \end{description}
\end{lemma}

\begin{lemma} \cite{W.J. Ning} \label{lem2}
Let $G$ be a connected graph and $G^{\prime}$ be a proper subgraph of $G$. Then $\rho(G^{\prime})<\rho(G)$.
\end{lemma}

The following lemma, due to Wu, Xiao and Hong \cite{B.F. Wu}, shows an edge transformation which increases the spectral radius.
\begin{lemma} \cite{B.F. Wu} \label{lem3}
Let $G$ be a connected graph and $\rho(G)$ be the spectral radius of $A(G)$. Let $u,v$ be two vertices of $G$. Suppose that $v_1,v_2, \ldots,v_s \in N_G(v) \backslash N_G(u)$ with $1\leq s\leq d_G(v)$, and $G^{\ast}$ is the graph obtained from $G$ by deleting the edges $v v_i$ and adding the edges $u v_i$ for $1 \leq i \leq s$. Let $\bm{x}$ be the Perron vector of $A(G)$. If $x(u) \geq x(v)$, then $\rho(G)<\rho(G^{\ast})$.
\end{lemma}

\begin{lemma}\label{lem4}\cite{D.D. Fan1}
Let $n=\sum_{i=1}^{t}n_i+s$. If $n_1\geq n_2\geq \cdots \geq n_t\geq p$ and $n_1<n-s-p(t-1)$, then
$$\rho(K_s\vee(K_{n_1}\cup K_{n_2} \cup\cdots \cup K_{n_t}))<\rho(K_s\vee(K_{n-s-p(t-1)}\cup (t-1)K_p)).$$
\end{lemma}

For any vertex $v\in V(G)$ and any subset $F\subseteq V(G)$, let $N_F(v)=N_G(v)\cap F$ and $d_F(v)=|N_G(v)\cap F|$. Now we shall give a proof of Theorem \ref{thm1}.

\vspace*{3mm}
\noindent{\textbf{Proof of Theorem \ref{thm1}}.}
Suppose that $G$ is a connected graph attaining the maximum spectral radius among all connected graph with minimum degree $\delta$ and $g$-good $r$-component connectivity $c\kappa_{g,r}(G)=k$.  By the definition of $g$-good $r$-component connectivity, there exists some nonempty subset $F\subseteq V(G)$ with $|F|=k$ such that $G-F$ is disconnected and there are at least $r$ components and each vertex $v\in V(G)\backslash F$ has at least $g$ neighbors. Let $B_1,B_2,\ldots,B_q$ be the components of $G-F$, where $q\geq r$, and let $|V(B_i)|=n_i$ for $1\leq i\leq q$. It is easy to see that $n_i\geq g+1$ for $1\leq i\leq q$. In fact, we can deduce that $q=r$ by Lemma \ref{lem2}. Now, we divide the proof into the following two cases.
\begin{case}
$k>\delta$. Choose a vertex $u\in V(G)$ such that $d_G(u)=\delta$. Suppose that $d_F(u)=t$. Next, we will discuss the case in two subcases.
\end{case}
\begin{subcase}\label{subcase1.1}
$\delta<g$.
\end{subcase}
By the definition of $g$-good $r$-component connectivity, we have $d_G(v)\geq g$ for $v\in V(G)\backslash F$. So $u\in F$.
We first assert that $t=\delta$. Otherwise, $t\leq\delta-1$. By the maximality of $\rho(G)$, we can obtain that
$G-u\cong K_{k-1}\vee(K_{n_1}\cup K_{n_2}\cup \cdots \cup K_{n_r})$. Notice that $d_F(u)=t$ and $|N_G(u)\cap (V(G)\backslash F)|\geq 1$.
Let $N_F(u)=\{w_1,w_2,\ldots,w_t\}$ and $F\backslash N_F[u]=\{w_{t+1},w_{t+2},\ldots,w_{k-1}\}$, and let $P_i=N_G(u)\cap V(K_{n_i})$ and $|P_i|=p_i$ for $1\leq i\leq r$. Then $\sum_{i=1}^{r}p_i=\delta-t$. Let $\bm{x}$ be the Perron vector of $A(G)$. By symmetry, we may say that $x(v)=x_i$ for $v\in P_i$ and $x(v)=x_i'$ for $v\in V(K_{n_i})\backslash P_i$, where $1\leq i\leq r$. Without loss of generality, assume that $x_1=\text{max}\{x_i \mid 1\leq i\leq r\}$. Note that $x(w_1)=x(w_i)$ for $2\leq i\leq t$ and $x(w_{t+1})=x(w_i)$ for $t+2\leq i\leq k-1$. Therefore, by $A(G)\bm{x}=\rho(G)\bm{x}$, we have
\begin{equation}\label{equ1}
\begin{aligned}
\rho(G)x(u)=tx(w_1)+\sum_{i=1}^rp_ix_i,
\end{aligned}
\end{equation}
\begin{equation}\label{equ2}
\begin{aligned}
\rho(G)x_i=tx(w_1)+(k-t-1) x(w_{t+1})+(p_i-1)x_i+(n_i-p_i)x_i^{\prime}+x(u),
\end{aligned}
\end{equation}
\begin{equation}\label{equ3}
\begin{aligned}
\rho(G)x_i^{\prime}=tx(w_1)+(k-t-1) x(w_{t+1})+p_ix_i+(n_i-p_i-1)x_i^{\prime},
\end{aligned}
\end{equation}
\begin{equation}\label{equ4}
\begin{aligned}
\rho(G)x(w_{t+1})=tx(w_1)+(k-t-2)x(w_{t+1})+\sum_{i=1}^{r} p_i x_i+\sum_{i=1}^{r}(n_i-p_i) x_i^{\prime},
\end{aligned}
\end{equation}
where $1\leq i\leq r$. From (\ref{equ1})-(\ref{equ3}), we get
\begin{equation*}
\begin{aligned}
& \rho(G)\left(\sum_{i=2}^{r}p_ix_i+\sum_{i=2}^{r}(n_i-p_i)x_i'-x(u)\right)\\
= &\left(\sum_{i=2}^{r}n_i-1\right)tx(w_1)+\sum_{i=2}^{r}n_i(k-t-1)x(w_{t+1})+\sum_{i=2}^{r}(n_i-1)p_ix_i+\sum_{i=2}^{r}(n_i-1)(n_i-p_i)x_i^{\prime}\\
&+\sum_{i=2}^{r}p_ix(u)-\sum_{i=1}^{r}p_ix_i \\
> & \sum_{i=2}^{r}n_i(k-t-1)x(w_{t+1})-\sum_{i=1}^{r}p_ix_i ~(\text{since}~n_i\geq g+1,~g>\delta,~t\geq 0~\text{and}~0\leq p_i<n_i~\text{for}~2\leq i\leq r) \\
\geq & (r-1)(g+1)(k-t-1)x(w_{t+1})-(\delta-t)x_1 \\
&(\text{since}~n_i\geq g+1,~\sum_{j=1}^rp_j=\delta-t~\text{and}~x_1\geq x_i~\text{for}~2\leq i\leq r) \\
> & (\delta-t)(x(w_{t+1})-x_1)~~(\text{since}~r\geq 2,~g>\delta~\text{and}~k\geq \delta+1).
\end{aligned}
\end{equation*}
Combining this with (\ref{equ2}) and (\ref{equ4}), we have
$$(\rho(G)+1)(x(w_{t+1})-x_{1})=\sum_{i=2}^{r}p_ix_i+\sum_{i=2}^{r}(n_i-p_i)x_i'-x(u) >\frac{(\delta-t)(x(w_{t+1})-x_1)}{\rho(G)},$$
from which we obtain
$$\frac{(\rho^{2}(G)+\rho(G)-\delta+t)(x(w_{t+1})-x_1)}{\rho(G)} > 0.$$
Note that $G$ contains $K_{k-1+n_i}$ as proper subgraphs, where $1\leq i\leq r$. Thus, $\rho(G)>\rho(K_{k-1+n_i})=k+n_i-2 >\delta$ due to $n_i\geq g+1>\delta+1$ and $k\geq \delta+1$. Combining this with $t\geq0$, we can deduce that $\rho^{2}(G)+\rho(G)-\delta+t>0$. This demonstrates that $x(w_{t+1})>x_1$.
Let $G'=G-\{uv \mid v\in P_i,1\leq i\leq r\}+\{uw_j \mid t+1\leq j\leq\delta\}$. Obviously, $G'\in\mathcal{G}_n^{k,\delta}$. According to Lemma \ref{lem3}, we have
$$\rho(G')>\rho(G),$$
which contradicts the maximality of $\rho(G)$. This implies that $t=\delta$.
In what follows, we shall prove that $G\cong G_{n,(g+1)^{r-1}}^{\delta,0}$. Without loss of generality, suppose that $n_1\geq n_2\geq \cdots\geq n_r$. By the maximality of $\rho(G)$ and $N_G[u]\subseteq F$, we have $G-u\cong K_{k-1}\vee(K_{n_1}\cup K_{n_2}\cup\cdots \cup K_{n_r})$, where $\sum_{i=1}^rn_i=n-k$ and $n_r\geq g+1$.
Denote $V(K_{n_i})=\{v_i^1,v_i^2,\ldots,v_i^{n_i}\}$. If $n_i=g+1$ for all $2\leq i\leq r$, then the result follows. If there exists $n_j\geq g+2$ for some $2\leq j\leq r$. Let $x(v_1^l)=x_1$ for $1\leq l\leq n_1$ and $x(v_j^i)=x_j$ for $2\leq j\leq r$ and $1\leq i\leq n_j$, then
$$(\rho(G)-n_j+1)(x_1-x_j)=(n_1-n_j)x_1\geq 0$$
due to $n_1\geq n_j$. Since $G$ contains $K_{k-1+n_j}$ as a proper subgraph and $k\geq \delta+1\geq 2$, it follows that $\rho(G)>\rho(K_{k-1+n_j})=k+n_j-2\geq n_j$, and hence $x_1\geq x_j$. Let $G''=G-\{v_j^iv_j^h\mid 1\leq i\leq g+1,~g+2\leq h\leq n_j\}+\{v_1^lv_j^h\mid 1\leq l\leq n_1,~g+2\leq h\leq n_j\}$. Then $G''-F$ contains $r$ components $K_{n-(r-1)(g+1)-k},K_{g+1},\ldots,K_{g+1}$, $|N_{G''}(v)\cap (V(G'')\backslash F)|\geq g$ for $v\in V(G'')\backslash F$ and $d_{G''}(u)=\delta$. So $G''\in\mathcal{G}_n^{k,\delta}$. According to Lemma \ref{lem3}, we have $\rho(G'')>\rho(G)$, which also leads to a contradiction. This indicates that $G\cong G_{n,(g+1)^{r-1}}^{\delta,0}$, as desired.

\vspace*{3mm}
\begin{subcase}
$\delta\geq g$. Recall that $u\in V(G)$ such that $d_G(u)=\delta$ and $d_F(u)=t$. We will prove the following four claims.
\end{subcase}
\begin{claim}\label{claim1}
$u\notin F$.
\end{claim}
Otherwise, $u\in F$. Using the similar analysis as the proof of subcase \ref{subcase1.1}, we can obtain that if $u\in F$, then $t=\delta$, i.e., $N_G(u)\subseteq F$. By the maximality of $\rho(G)$, we can deduce that
$G-u\cong K_{k-1}\vee(K_{n_1}\cup K_{n_2}\cup\cdots \cup K_{n_r})$, where $\sum_{i=1}^rn_i=n-k$ and $n_i\geq g+1$. Without loss of generality, we may consider that $n_1\geq n_2\geq \cdots\geq n_r$.
Let $N_F(u)=\{w_1,w_2,\ldots,w_\delta\}$, $F\backslash N_F[u]=\{w_{\delta+1},w_{\delta+2},\ldots,w_{k-1}\}$
and $V(K_{n_i})=\{v_i^1, v_i^2,\ldots, v_i^{n_i}\}$ for $1\leq i\leq r$. Assume that $\bm{x}$ is the Perron vector of $A(G)$. By symmetry, we say $x(v)=x_i$ for $v\in V(K_{n_i})$, where $1\leq i\leq r$. Observe that $x(w_1)=x(w_i)$ for $2\leq i\leq \delta$ and $x(w_{\delta+1})=x(w_i)$ for $\delta+2\leq i\leq k-1$.
Since $N_G(v_1^i)\backslash \{w_1\}\subset N_G(w_1)\backslash \{v_1^i\}$ for $1\leq i\leq n_1$, it follows that $x_1<x(w_1)$ by Lemma \ref{lem1}. Furthermore, if $n_i\geq n_j$, then we can deduce that $(\rho(G)-n_j+1)(x_i-x_j)=(n_i-n_j)x_i\geq 0$. Note that $G$ contains $K_{k-1+n_j}$ as a proper subgraph and $k\geq \delta+1\geq 2$, then $\rho(G)>\rho(K_{k-1+n_j})=k+n_j-2\geq n_j$, and hence $x_i\geq x_j$. Moreover,
\begin{equation*}
\begin{aligned}
\rho(G)x(u)=\delta x(w_1),
\end{aligned}
\end{equation*}
\begin{equation*}
\begin{aligned}
\rho(G)x_i=\delta x(w_1)+(k-\delta-1)x(w_{\delta+1})+(n_i-1)x_i,
\end{aligned}
\end{equation*}
\begin{equation*}
\begin{aligned}
\rho(G)x(w_1)=(\delta-1)x(w_1)+(k-\delta-1)x(w_{\delta+1})+x(u)+\sum_{i=1}^rn_ix_i,
\end{aligned}
\end{equation*}
where $1\leq i\leq r$, from which we have
\begin{align} \label{equ5}
 &(n-k-n_r)x_r+gx(u)-gx(w_1)\nonumber\\
=&\frac{(\rho(G)+1)((n-k-n_r)x_r+gx(u)-gx(w_1))}{\rho(G)+1}\nonumber\\
=&\frac{(n-k-n_r)\delta x(w_1)+(n-k-n_r-g)(k-\delta-1)x(w_{\delta+1})+(n-k-n_r)n_rx_r-g\sum_{i=1}^{r}n_ix_i}{\rho(G)+1}\nonumber\\
\geq&\frac{(n-k-n_r)\delta x(w_1)+(n-k-n_r)n_rx_r-gn_rx_r-g\sum_{i=1}^{r-1}n_ix_i}{\rho(G)+1} ~~(\text{since}~n-k-n_r>g~\text{and}~k\geq \delta+1)\nonumber\\
\geq&\frac{(n-k-n_r)\delta x(w_1)+(n-k-n_r-g)n_rx_r-g(n-k-n_r)x_1}{\rho(G)+1} \nonumber\\ &(\text{since}~\sum_{i=1}^{r-1}n_i=n-k-n_r~\text{and}~x_1\geq x_j~\text{for}~2\leq j\leq r-1)\nonumber\\
=&\frac{(n-k-n_r)(\delta x(w_1)-gx_1)+(n-k-n_r-g)n_rx_r}{\rho(G)+1} \nonumber\\
>&0~~(\text{since}~\delta\geq g,~x(w_1)>x_1~\text{and}~n-k-n_r>g)
\end{align}
and
\begin{equation}\label{equ6}
\begin{aligned}
 x_{r-1}-x(u)=\frac{(k-\delta-1)x(w_{\delta+1})+(n_{r-1}-1)x_{r-1}}{\rho(G)}\geq 0.
\end{aligned}
\end{equation}
Suppose that $E_1=\{v_r^1v\mid v\in V(G)\backslash (F\cup V(K_{n_r}))\}+\{uv_r^i\mid 2\leq i\leq g+1\}$ and $E_2=\{uw_j\mid \delta-g+1\leq j\leq\delta\}$. Let $G_1=G+E_1-E_2$ and $F_1=F-\{u\}+\{v_r^1\}$. Then $|F_1|=k$, $G_1-F_1$ contains $r$ components $K_{n_1},K_{n_2},\ldots,K_{n_{r-1}},K_{n_r-1}\cup \{u\}$ and $d_{G_1}(u)=\delta$. One can verify that $G_1\in\mathcal{G}_n^{k,\delta}$.
Let $\bm{y}$ be the Perron vector of $A(G_1)$. By symmetry, $y(v)=y_i$ for $v\in V(K_{n_i})$, where $1\leq i\leq r-1$, $y(v_r^i)=y_r$ for $2\leq i\leq g+1$ and $y(v_r^j)=y_r'$ for $g+2\leq j\leq n_r$. Notice that $y(w_1)=y(w_i)$ for $2\leq i\leq \delta-g$ and $y(v_r^1)=y(w_i)$ for $\delta-g+1\leq i\leq k-1$. According to Lemma \ref{lem1}, we get $y(u)<y_r$ and $y_r'<y_r$, then
\begin{equation*}
\begin{aligned}
\rho(G_1)y_i=(\delta-g)y(w_1)+(k-\delta+g)y(v_r^1)+(n_i-1)y_i,
\end{aligned}
\end{equation*}
\begin{equation*}
\begin{aligned}
\rho(G_1)y_r=&(\delta-g)y(w_1)+(k-\delta+g)y(v_r^1)+(g-1)y_r+(n_r-g-1)y_r'+y(u)\\
  <&(\delta-g)y(w_1)+(k-\delta+g)y(v_r^1)+(n_r-1)y_r,
\end{aligned}
\end{equation*}
where $1\leq i\leq r-1$, from which we have
\begin{equation*}
\begin{aligned}
(\rho(G_1)-n_r+1)(y_i-y_r)>(n_i-n_r)y_i\geq 0,
\end{aligned}
\end{equation*}
due to $n_i\geq n_r$. Observe that $G_1$ contains $K_{k+n_r-1}$ as a proper subgraph, $\rho(G_1)>\rho(K_{k+n_r-1})=k+n_r-2\geq n_r$. Thus $y_i>y_r$ for $1\leq i\leq r-1$. Similarly, if $n_i\geq n_j$, then $y_i\geq y_j$. Combining these with (\ref{equ5}) and (\ref{equ6}), we obtain
\begin{equation*}
\begin{aligned}
&\bm{y}^T(\rho(G_1)-\rho(G))\bm{x}\\
=& \bm{y}^T(A(G_1)-A(G))\bm{x}\\
=& \sum_{v_r^1v\in E_1}(x(v_r^1)y(v)+x(v)y(v_r^1))+\sum_{uv_r^i\in E_1}(x(u)y(v_r^i)+x(v_r^i)y(u))-\sum_{uw_j\in E_2}(x(u)y(w_j)+x(w_j)y(u))\\
\geq & (n-k-n_r)(x_ry_{r-1}+x_{r-1}y(v_r^1))+g(x(u)y_r+x_ry(u))-g(x(u)y(v_r^1)+x(w_1)y(u))\\
&(\text{since}~x_i\geq x_{r-1}~\text{and}~y_i\geq y_{r-1}~\text{for}~1\leq i\leq r-2)\\
> & (n-k-n_r)x_ry_{r-1}+gx(u)y_r-gx(w_1)y(u)~~(\text{since}~n-k-n_r>g,~x_{r-1}\geq x(u)~\text{and}~g\geq 0)\\
>& y(u)((n-k-n_r)x_r+gx(u)-gx(w_1))~~(\text{since}~y_{r-1}>y_r>y(u))\\
>& 0~~(\text{by}~(\ref{equ5})),
\end{aligned}
\end{equation*}
and hence $\rho(G_1)>\rho(G)$, which contradicts the maximality of $\rho(G)$. This demonstrates that $u\notin F$, proving Claim \ref{claim1}.

\vspace{2mm}
By Claim \ref{claim1}, $u\notin F$, and so $u\in V(B_i)$ for some $1\leq i\leq r$. Without loss of generality, we may assume that $u\in V(B_r)$. Then $d_{B_r}(u)\geq g$ and $t=d_{F}(u)=\delta-d_{B_r}(u)\leq\delta-g$. We can get the following claim.

\begin{claim}\label{claim2}
$t=\delta-g$.
\end{claim}
If $\delta=g$, then $t=0$. Clearly, the result holds. Next, we consider $\delta>g$. Suppose on the contrary that $t< \delta-g$, i.e., $d_{B_r}(u)=\delta-t\geq g+1$. Let $N_F(u)=\{w_1,w_2,\ldots, w_t\}$, $F\backslash N_F(u)=\{w_{t+1},w_{t+2},\ldots, w_k\}$ and $N_{B_r}(u)=\{v_r^1,v_r^2,\ldots,v_r^{\delta-t}\}$. Again by the maximality of $\rho(G)$, we can deduce that $G-u\cong K_k\vee(K_{n_1}\cup\cdots \cup K_{n_{r-1}} \cup K_{n_r-1})$. Assume that $\bm{x}$ is the Perron vector of $A(G)$. By symmetry, we may say $x(v)=x_r$ for $v \in N_{B_r}(u)$, $x(v)=x_{r}^{\prime}$ for $v \in V(B_r)\backslash N_{B_r}[u]$ and $x(v)=x_i$ for $v \in V(K_{n_i})$, where $1\leq i\leq r-1$. Notice that $x(w_1)=x(w_i)$ for $2\leq i\leq t$ and $x(w_{t+1})=x(w_i)$ for $t+2 \leq i \leq k$, then
\begin{equation*}
\begin{aligned}
\rho(G)x(w_{t+1})=tx(w_1)+(k-t-1)x(w_{t+1})+(\delta-t)x_r+(n_r-1-(\delta-t)) x_r^{\prime}+\sum_{i=1}^{r-1}n_ix_i,
\end{aligned}
\end{equation*}
\begin{equation*}
\begin{aligned}
\rho(G)x_r=tx(w_1)+(k-t)x(w_{t+1})+(\delta-t-1)x_r+(n_r-1-(\delta-t))x_r^{\prime}+x(u),
\end{aligned}
\end{equation*}
from which we obtain that
\begin{align*}
&(\rho(G)+1)(x(w_{t+1})-x_r)\\
=& \sum_{i=1}^{r-1}n_ix_i-x(u) \\
=& \frac{1}{\rho(G)}\left(\sum_{i=1}^{r-1}n_i(\rho(G)x_i)-\rho(G)x(u)\right)\\
= & \frac{\sum_{i=1}^{r-1}n_i(tx(w_1)+(k-t)x(w_{t+1})+(n_i-1)x_i)-(tx(w_1)+(\delta-t)x_r)}{\rho(G)} \\
= & \frac{(\sum_{i=1}^{r-1}n_i-1)tx(w_1)+\sum_{i=1}^{r-1}n_i(k-t)x(w_{t+1})+\sum_{i=1}^{r-1}n_i(n_i-1)x_i-(\delta-t)x_r}{\rho(G)}\\
\geq & \frac{(r-1)(g+1)(k-t)x(w_{t+1})-(\delta-t)x_r}{\rho(G)}~~ (\text{since}~n_i\geq g+1~\text{for}~1\leq i\leq r-1,~t\geq 0~\text{and}~g\geq 0) \\
>&\frac{(\delta-t)(x(w_{t+1})-x_r)}{\rho(G)}~~(\text{since}~r\geq 2,~g\geq 0~\text{and}~k>\delta).
\end{align*}
It follows that $\frac{(\rho^{2}(G)+\rho(G)-\delta+t)(x(w_{t+1})-x_r)}{\rho(G)}>0$. Note that $G$ contains $K_k$ as a proper subgraph. Thus, $\rho(G)>\rho(K_k)=k-1\geq\delta$. Combining this with $t\geq 0$, we can deduce that $\rho^{2}(G)+\rho(G)-\delta+t>0$. This suggests that $x(w_{t+1})>x_r$. Let $G_2=G-\{uv_r^i \mid g+1\leq i \leq\delta-t\}+\{uw_j\mid t+1\leq j\leq\delta-g\}$. It is obvious that $G_2\in \mathcal{G}_n^{k,\delta}$. According to Lemma \ref{lem3}, we have $\rho(G_2)>\rho(G)$, which contradicts the maximality of $\rho(G)$. Therefore, $t=\delta-g$, completing the proof of Claim \ref{claim2}.

\vspace{2mm}
Without loss of generality, suppose that $n_1\geq n_2\geq\cdots\geq n_r$. If $n_1=n_r$, then $u\in V(B_i)$ for $1\leq i\leq r$. If $n_1>n_r$, then we have the following claim.
\begin{claim}\label{claim3}
$u\in V(B_i)$ for $n_i=n_r$.
\end{claim}
If $g=0$, then $t=\delta$. Clearly, $u\in V(B_i)$ for $n_i=n_r=1$. Therefore, the result holds. Next, we consider $g\geq 1$. Suppose on the contrary that $u\in V(B_j)$ for $n_j>n_r$. Without loss of generality, assume that $u\in V(B_1)$. We obtain $G-u\cong K_k\vee(K_{n_1-1}\cup K_{n_2}\cup \cdots \cup K_{n_r})$ by the maximality of $\rho(G)$. Let $N_F(u)=\{w_1,w_2,\ldots,w_{\delta-g}\}$ and $F\backslash N_F(u)=\{w_{\delta-g+1},w_{\delta-g+2},\ldots,w_k\}$. Let $\bm{x}$ be the Perron vector of $A(G)$. By symmetry, $x(v)=x_1$ for $v \in N_{B_1}(u)$, $x(v)=x_{1}^{\prime}$ for $v \in V(B_1) \backslash N_{B_1}[u]$ and $x(v)=x_i$ for $v \in V(K_{n_i})$, where $2\leq i\leq r$. Notice that $x(w_1)=x(w_i)$ for $2\leq i\leq \delta-g$ and $x(w_{\delta-g+1})=x(w_i)$ for $\delta-g+2 \leq i \leq k$, then
$$\rho(G)x_r=(\delta-g)x(w_1)+(k-\delta+g)x(w_{\delta-g+1})+(n_r-1)x_r,$$
$$\rho(G)x_1'=(\delta-g)x(w_1)+(k-\delta+g)x(w_{\delta-g+1})+gx_1+(n_1-g-2)x_1',$$
$$\rho(G)x_1=(\delta-g)x(w_1)+(k-\delta+g)x(w_{\delta-g+1})+(g-1)x_1+(n_1-g-1)x_1'+x(u),$$
from which we have
\begin{equation}\label{equ7}
\begin{aligned}
x_1=x_1'+\frac{x(u)}{\rho(G)+1}.
\end{aligned}
\end{equation}
Since $G$ contains $K_{k+n_1-1}$ as a proper subgraph and $n_1>n_r$, it follows that $\rho(G)>\rho(K_{k+n_1-1})=k+n_1-2$ and $\rho(G)-n_r+1\geq\rho(G)-n_1+2>0$. Then we get
\begin{equation}\label{equ8}
\begin{aligned}
x_1'=\frac{\rho(G)-n_r+1}{\rho(G)-n_1+2}x_r+\frac{g(x_1-x_1')}{\rho(G)-n_1+2}> x_r.
\end{aligned}
\end{equation}
Moreover, we have $x_1>x(u)$ by Lemma \ref{lem1} and
\begin{equation*}
\begin{aligned}
\rho(G)x(u)=(\delta-g)x(w_1)+gx_1,
\end{aligned}
\end{equation*}
\begin{equation*}
\begin{aligned}
\rho(G)x(w_{\delta-g+1})=(\delta-g)x(w_1)+(k-\delta+g-1)x(w_{\delta-g+1})+gx_1+(n_1-g-1)x_1'+\sum_{i=2}^rn_ix_i,
\end{aligned}
\end{equation*}
from which we get that
\begin{align} \label{equ9}
x_r-x(u)=& \frac{1}{\rho(G)}(\rho(G)x_r-\rho(G)x(u))\nonumber\\
=& \frac{(k-\delta+g)x(w_{\delta-g+1})+(n_r-1)x_r-gx_1}{\rho(G)}\nonumber\\
>& \frac{g(x(w_{\delta-g+1})+x_r-x_1)}{\rho(G)}~~(\text{since}~k>\delta~\text{and}~n_r\geq g+1)\nonumber\\
=& \frac{g}{\rho^2(G)}(\rho(G)x(w_{\delta-g+1})+\rho(G)x_r-\rho(G)x_1)\nonumber\\
=& \frac{g}{\rho^2(G)}[(\delta-g)x(w_1)+(k-\delta+g-1)x(w_{\delta-g+1})+(n_r-1)x_r+\sum_{i=2}^rn_ix_i+x_1-x(u)]\nonumber\\
>& 0~~(\text{since}~\delta\geq g,~k\geq \delta+1,~g\geq 1,~n_r\geq g+1~\text{and}~x_1>x(u)).
\end{align}
Denote $V(B_1)=\{v_1^1,v_1^2,\ldots,v_1^{n_1-1},u\}$, $N_{B_1}(u)=\{v_1^1,v_1^2,\ldots,v_1^g\}$ and $V(K_{n_i})=\{v_i^1,v_i^2,\ldots,v_i^{n_i}\}$ for $2\leq i\leq r$.
Suppose that $E_1'=\{uv_r^i\mid 1\leq i\leq g\}+\{v_r^{n_r}v_1^j\mid 1\leq j\leq n_1-1\}$ and $E_2'=\{uv_1^i\mid 1\leq i\leq g\}+\{v_r^{n_r}v_r^j\mid 1\leq j\leq n_r-1\}$. Let $G_3=G+E_1'-E_2'$, then $V(B_1')=\{v_1^1,\ldots,v_1^{n_1-1}, v_r^{n_r}\}$ and $V(B_r')=\{v_r^1,\ldots,v_r^{n_r-1},u\}$. So $G_3-F$ contains $r$ components $K_{n_1},K_{n_2},\ldots,K_{n_{r-1}},K_{n_r-1}\cup \{u\}$ and $d_{G_3}(u)=\delta$. Obviously, $G_3\in\mathcal{G}_n^{k,\delta}$. Combining this with (\ref{equ7})-(\ref{equ9}), we have
\begin{equation*}
\begin{aligned}
&\rho(G_3)-\rho(G)\\
\geq & \bm{x}^T(A(G_3)-A(G))\bm{x}\\
= & 2\sum_{uv_r^i\in E_1'}x(u)x(v_r^i)+2\sum_{v_r^{n_r}v_1^j\in E_1'}x(v_r^{n_r})x(v_1^j)-2\sum_{uv_1^i \in E_2'}x(u)x(v_1^i)-2\sum_{v_r^{n_r}v_r^j\in E_2'}x(v_r^{n_r})x(v_r^j)\\
= & 2gx(u)x_r+2gx_rx_1+2(n_1-g-1)x_rx_1'-2gx(u)x_1-2(n_r-1)x_rx_r\\
\geq & 2gx(u)x_r+\frac{2gx_rx(u)}{\rho(G)+1}+2n_rx_rx_1'-2gx(u)x_1'-\frac{2gx(u)x(u)}{\rho(G)+1}-2(n_r-1)x_rx_r
~~(\text{since}~n_1\geq n_r+1)\\
> & 2gx(u)x_r+2n_rx_rx_1'-2gx(u)x_1'-2n_rx_rx_r+2x_rx_r~~(\text{since}~x_r>x(u)~\text{and}~g\geq 1)\\
> & 2n_rx_r(x_1'-x_r)-2gx(u)(x_1'-x_r)\\
=& 2(n_rx_r-gx(u))(x_1'-x_r)\\
>& 0~~(\text{since}~n_r\geq g+1,~x_r>x(u)~\text{and}~x_1'> x_r).
\end{aligned}
\end{equation*}
Hence $\rho(G_3)>\rho(G)$, which contradicts the maximality of $\rho(G)$. This implies that $u\in V(B_i)$ for $n_i=n_r$, proving Claim \ref{claim3}.

\vspace{2mm}
In what follows, we shall prove that $G\cong G_{n,(g+1)^{r-1}}^{\delta-g,g}$. In fact, by the maximality of $\rho(G)$ and Claims \ref{claim1}-\ref{claim3}, we get $G-u\cong K_k\vee(K_{n_1}\cup K_{n_2}\cup \cdots \cup K_{n_{r-1}}\cup K_{n_r-1})$ for $n_1\geq n_2\geq \cdots\geq n_r\geq g+1$, where $\sum_{i=1}^rn_i=n-k$. Then we have the following claim.
\begin{claim} \label{claim4}
$n_i=g+1$ for $2\leq i\leq r$.
\end{claim}
Otherwise, there exists $n_j\geq g+2$ for some $2\leq j\leq r$. Let $V(K_{n_i})=\{v_i^1,v_i^2,\ldots,v_i^{n_i}\}$ for $1\leq i\leq r-1$, $N_{B_r}(u)=\{v_r^1,v_r^2,\ldots,v_r^g\}$ and $V(B_r)\backslash N_{B_r}[u]=\{v_r^{g+1},v_r^{g+2},\ldots,v_r^{n_r-1}\}$. Assume that $\bm{x}$ is the Perron vector of $A(G)$. By symmetry,
$x(v_i^l)=x_i$ for $1\leq i\leq r-1$ and $1\leq l\leq n_i$, $x(v_r^i)=x_r$ for $1\leq i\leq g$ and $x(v_r^j)=x_r'$ for $g+1\leq j\leq n_r-1$. Observe that $x(w_1)=x(w_i)$ for $2\leq i\leq \delta-g$ and $x(w_{\delta-g+1})=x(w_i)$ for $\delta-g+2 \leq i \leq k$. According to Lemma \ref{lem1}, we have $x_r'<x_r$ and $x(u)<x_r$. Furthermore,
$$\rho(G)x_1=(\delta-g)x(w_1)+(k-\delta+g)x(w_{\delta-g+1})+(n_1-1)x_1,$$
\begin{equation*}
\begin{aligned}
\rho(G)x_r&=(\delta-g)x(w_1)+(k-\delta+g)x(w_{\delta-g+1})+(g-1)x_r+(n_r-g-1)x_r'+x(u)\\
           &<(\delta-g)x(w_1)+(k-\delta+g)x(w_{\delta-g+1})+(n_r-1)x_r,
\end{aligned}
\end{equation*}
from which we obtain that
$$(\rho(G)-n_r+1)(x_1-x_r)>(n_1-n_r)x_1\geq 0$$
due to $n_1\geq n_r$. Note that $G$ contains $K_{k+n_1}$ as a proper subgraph, it follows that $\rho(G)>\rho(K_{k+n_1})=k+n_1-1>n_r$ due to $k\geq \delta+1\geq 2$. Thus $x_1> x_r$. Similarly, we have $x_1\geq x_j$ for $2\leq j\leq r-1$. Let $G_4=G-\{v_j^iv_j^l \mid 2\leq j\leq r-1,1\leq i\leq g+1,g+2\leq l\leq n_j\}+\{v_1^pv_j^l \mid 1\leq p\leq n_1,2\leq j\leq r-1,g+2\leq l\leq n_j\}-\{v_r^iv_r^j\mid 1\leq i\leq g,g+1\leq j\leq n_r-1\}+\{v_1^pv_r^j\mid 1\leq p\leq n_1,g+1\leq j\leq n_r-1\}$. One can verify that $G_4\in\mathcal{G}_n^{k,\delta}$. According to Lemma \ref{lem3}, we have $\rho(G_4)>\rho(G)$, which also leads to a contradiction. \par
Based on the above results, we can conclude that $G\cong G_{n,(g+1)^{r-1}}^{\delta-g,g}$, as desired.

\vspace*{3mm}
\begin{case}
$k\leq \delta$. Recall that $u\in V(G)$ such that $d_G(u)=\delta$ and $d_F(u)=t$. Next, we will discuss the case in three subcases.
\end{case}
\begin{subcase} \label{sub2.1}
$\delta<g$.
\end{subcase}
By the definition of $g$-good $r$-component connectivity, we have $d_G(v)\geq g$ for $v\in V(G)\backslash F$. So $u\in F$. Then we consider the following two ways according to the value of $k$. \par
\textbf{In the case of $k\geq 2$.}
 Since $u\in F$ and $|F|=k\leq \delta$, we have $t\leq k-1$, then we will get the following claim.
\begin{claim} \label{claim5}
$t=k-1$.
\end{claim}
Otherwise, $t\leq k-2$. Then $|N_G(u)\cap (V(G)\backslash F)|\geq 2$. By the maximality of $\rho(G)$, we can obtain that
$G-u\cong K_{k-1}\vee(K_{n_1}\cup K_{n_2}\cup \cdots \cup K_{n_r})$.
Let $N_F(u)=\{w_1,w_2,\ldots,w_t\}$ and $F\backslash N_F[u]=\{w_{t+1},w_{t+2},\ldots,w_{k-1}\}$, and let $P_i=N_G(u)\cap V(K_{n_i})$ and $|P_i|=p_i$ for $1\leq i\leq r$. Then $\sum_{i=1}^rp_i=\delta-t$. Let $\bm{x}$ be the Perron vector of $A(G)$. By symmetry, we may say $x(v)=x_i$ for $v\in P_i$ and $x(v)=x_i'$ for $v\in V(K_{n_i})\backslash P_i$, where $1\leq i\leq r$. Without loss of generality, assume that $x_1=\text{max}\{x_i\mid 1\leq i\leq r\}$. Note that $x(w_1)=x(w_i)$ for $2\leq i\leq t$ and $x(w_{t+1})=x(w_i)$ for $t+2\leq i\leq k-1$. Using the similar analysis as the proof of subcase \ref{subcase1.1}, we get $x(w_{t+1})>x_1$.
Let $P_1\cup P_2\cup\cdots \cup P_r=\{v_1,\ldots,v_{k-t-1},\ldots,v_{\delta-t}\}$
and $G^{\star}=G-\{uv_i\mid 1\leq i\leq k-t-1\}+\{uw_j\mid t+1\leq j\leq k-1\}$. Clearly, $G^{\star}\in\mathcal{G}_n^{k,\delta}$. According to Lemma \ref{lem3}, we have
$$\rho(G^{\star})>\rho(G),$$
which contradicts the maximality of $\rho(G)$. This implies that $t=k-1$.

~\par
By Claim \ref{claim5}, we get $|N_G(u)\cap (V(G)\backslash F)|=\delta-k+1$, i.e., $\sum_{i=1}^rp_i=\delta-k+1$. We still consider that $x_1=\text{max}\{x_i\mid 1\leq i\leq r\}$. Then we have the following claim.
\begin{claim} \label{claim6}
$p_1=\delta-k+1$.
\end{claim}
Otherwise, $p_1\leq\delta-k$. Then $\sum_{i=2}^rp_i\geq 1$. Denote $P_i=\{v_i^1,v_i^2,\ldots, v_i^{p_i}\}$ and $V(B_i)\backslash P_i=\{v_i^{p_i+1},v_i^{p_i+2},\ldots, v_i^{n_i}\}$, where $1\leq i\leq r$. By the maximality of $\rho(G)$, we have
$G-u\cong K_{k-1}\vee(K_{n_1}\cup K_{n_2}\cup \cdots \cup K_{n_r})$ and $d_F(u)=k-1$. Let $\bm{x}$ be the Perron vector of $A(G)$. By symmetry, $x(w_1)=x(w_j)$ for $2\leq j\leq k-1$, $x(v)=x_i$ for $v\in P_i$ and $x(v)=x_i'$ for $v\in V(K_{n_i})\backslash P_i$, where $1\leq i\leq r$. We first assert that $n_1\geq n_i$ for $2\leq i\leq r$. If not, there exists $n_j> n_1$ for some $2\leq j\leq r$.
By $A(G)\bm{x}=\rho(G)\bm{x}$, we obtain
\begin{equation*}
\begin{aligned}
\rho(G) x_i=(k-1)x(w_1)+(p_i-1)x_i+(n_i-p_i)x_i'+x(u),
\end{aligned}
\end{equation*}
\begin{equation*}
\begin{aligned}
\rho(G) x_i'=(k-1)x(w_1)+p_ix_i+(n_i-p_i-1)x_i',
\end{aligned}
\end{equation*}
where $1\leq i\leq r$, from which we have
\begin{equation} \label{equ10}
\begin{aligned}
x_i=x_i'+\frac{x(u)}{\rho(G)+1}
\end{aligned}
\end{equation}
and
\begin{equation} \label{equ11}
\begin{aligned}
(\rho(G)+1)(x_1-x_j)=p_1x_1+(n_1-p_1)x_1'-p_jx_j-(n_j-p_j)x_j'\geq 0
\end{aligned}
\end{equation}
due to $x_1\geq x_j$. Let $G_1=G+\{v_j^{n_j}v \mid v\in V(K_{n_1})\}-\{v_j^{n_j}v' \mid v'\in V(K_{n_j})\backslash \{v_j^{n_j}\}\}$. Clearly, $G_1\in \mathcal{G}_n^{k,\delta}$. Combining this with (\ref{equ11}), we get
\begin{equation*}
\begin{aligned}
\rho(G_1)-\rho(G)&\geq \bm{x}^T(A(G_1)-A(G))\bm{x}\\
&=2x_j'p_1x_1+2x_j'(n_1-p_1)x_1'-2x_j'p_jx_j-2x_j'(n_j-p_j-1)x_j'\\
&=2x_j'[p_1x_1+(n_1-p_1)x_1'-p_jx_j-(n_j-p_j-1)x_j']\\
&>0~~(\text{by}~\ref{equ11}).
\end{aligned}
\end{equation*}
Hence $\rho(G_1)>\rho(G)$, which contradicts the maximality of $\rho(G)$. This suggests that $n_1\geq n_i$ for $2\leq i\leq r$.
Furthermore, let $G_2=G-\{uv\mid v\in P_i,~2\leq i\leq r\}+\{uv_1^j \mid p_1+1\leq j\leq \delta-k+1\}$. One can easily verify that $G_2\in \mathcal{G}_n^{k,\delta}$ and $p_1=\delta-k+1$. If $\rho(G_2)>\rho(G)$, then the results follows. If not, assume that $\bm{y}$ is the Perron vector of $A(G_2)$. By symmetry, we may say $y(v_j^i)=y_j$ for $2\leq j\leq r$ and $1\leq i\leq n_j$, $y(v_1^i)=y_1$ for $1\leq i\leq \delta-k+1$ and $y(v_1^j)=y_1'$ for $\delta-k+2\leq j\leq n_1$. Note that $y(w_1)=y(w_i)$ for $2\leq i\leq k-1$ and $y_1>y_1'$ by Lemma \ref{lem1}, then
\begin{equation*}
\begin{aligned}
\rho(G_2)y_1=(k-1)y(w_1)+y(u)+(\delta-k)y_1+(n_1-\delta+k-1)y_1',
\end{aligned}
\end{equation*}
\begin{equation*}
\begin{aligned}
\rho(G_2)y_1'=(k-1)y(w_1)+(\delta-k+1)y_1+(n_1-\delta+k-2)y_1'>(k-1)y(w_1)+(n_1-1)y_1',
\end{aligned}
\end{equation*}
\begin{equation*}
\begin{aligned}
\rho(G_2)y_j=(k-1)y(w_1)+(n_j-1)y_j,
\end{aligned}
\end{equation*}
where $2\leq j\leq r$, from which we get
\begin{equation} \label{equ12}
\begin{aligned}
y_1=y_1'+\frac{y(u)}{\rho(G_2)+1}.
\end{aligned}
\end{equation}
Recall that $n_1\geq n_j$ for $2\leq j\leq r$, we can deduce that
$(\rho(G_2)-n_1+1)(y_1'-y_j)>(n_1-n_j)y_j\geq 0$.
Observe that $G_2$ contains $K_{n_1}$ as a proper subgraph, thus $\rho(G_2)>\rho(K_{n_1})=n_1-1$. It follows that $y_1'>y_j$ for $2\leq j\leq r$. Combining this with (\ref{equ10}) and (\ref{equ12}), we obtain
\begin{equation*}
\begin{aligned}
&\bm{y}^T(\rho(G_2)-\rho(G))\bm{x}\\
=& \bm{y}^T(A(G_2)-A(G))\bm{x}\\
=& (\delta-k+1-p_1)(x(u)y(v_1^j)+x(v_1^j)y(u)-x(u)y(v)-x(v)y(u))\\
>& (\delta-k+1-p_1)(x(u)y_1+x_1'y(u)-x(u)y_1'-x_1y(u))~~(\text{since}~y_1'>y_j~\text{and}~x_1\geq x_j~\text{for}~2\leq j\leq r)\\
=& (\delta-k+1-p_1)\left(\frac{x(u)y(u)}{\rho(G_2)+1}-\frac{x(u)y(u)}{\rho(G)+1}\right)~~(\text{by}~(\ref{equ10}) ~\text{and}~(\ref{equ12}))\\
=& (\delta-k+1-p_1)\cdot \frac{(\rho(G)-\rho(G_2))x(u)y(u)}{(\rho(G)+1)(\rho(G_2)+1)}\\
\geq& 0~~(\text{since}~\rho(G_2)\leq \rho(G)),
\end{aligned}
\end{equation*}
a contradiction. This implies that $p_1=\delta-k+1$, proving Claim \ref{claim6}.

~\par
Without loss of generality, we may assume that $n_1\geq n_2\geq\cdots\geq n_r$. By the maximality of $\rho(G)$, we can obtain that
$G-u\cong K_{k-1}\vee(K_{n_1}\cup K_{n_2}\cup \cdots \cup K_{n_r})$, where $n_r\geq g+1$ and $\sum_{i=1}^rn_i=n-k$. Then we have the following claim.
\begin{claim}\label{claim7}
$n_i=g+1$ for $2\leq i\leq r$.
\end{claim}
Otherwise, there exists $n_j\geq g+2$ for some $2\leq j\leq r$. Denote $V(K_{n_i})=\{v_i^1,v_i^2,\ldots,v_i^{n_i}\}$ for $1\leq i\leq r$ and $P_1=\{v_1^1,v_1^2,\ldots,v_1^{\delta-k+1}\}$. Assume that $\bm{x}$ is the Perron vector of $A(G)$. By symmetry, we say
$x(v_j^i)=x_j$ for $2\leq j\leq r$ and $1\leq i\leq n_j$, $x(v_1^i)=x_1$ for $1\leq i\leq \delta-k+1$ and  $x(v_1^j)=x_1'$ for $\delta-k+2\leq j\leq n_1$. Notice that $x(w_1)=x(w_i)$ for $2\leq i\leq k-1$ and $x_1>x_1'$ by Lemma \ref{lem1}. Furthermore,
$$\rho(G)x_1'=(k-1)x(w_1)+(\delta-k+1)x_1+(n_1-\delta+k-2)x_1'>(k-1)x(w_1)+(n_1-1)x_1',$$
$$\rho(G)x_j=(k-1)x(w_1)+(n_j-1)x_j,$$
where $2\leq j\leq r$, from which we obtain that
$$(\rho(G)-n_1+1)(x_1'-x_j)>(n_1-n_j)x_j\geq 0$$
due to $n_1\geq n_j$. Note that $G$ contains $K_{n_1}$ as a proper subgraph, it follows that $\rho(G)>\rho(K_{n_1})=n_1-1$. Thus $x_1'> x_j$ and $x_1>x_j$. Let $G'=G-\{v_j^iv_j^l \mid 2\leq j\leq r,1\leq i\leq g+1,g+2\leq l\leq n_j\}+\{v_1^hv_j^l \mid 1\leq h\leq n_1,2\leq j\leq r,g+2\leq l\leq n_j\}$. It is not difficult to see that $G'\in\mathcal{G}_n^{k,\delta}$. According to Lemma \ref{lem3}, we have $\rho(G')>\rho(G)$, which contradicts the maximality of $\rho(G)$. \par
Based on the above results, we can conclude that $G\cong G_{n,(g+1)^{r-1}}^{k-1,\delta-k+1}$, as desired.

~\par
\textbf{In the case of $k=1$.} Recall that $u\in F$ and $q=r$. It is worth noting that $G$ is a connected graph, then $\delta\geq r$, $\sum_{i=1}^rp_i=\delta$ and $p_i\geq 1$ for $1\leq i\leq r$.
Without loss of generality, assume that $n_1\geq n_2\geq\cdots\geq n_r$. Notice that $n_r\geq g+1$.
By the maximality of $\rho(G)$ and using the similar analysis as the proof of Claims \ref{claim6} and \ref{claim7}, we can obtain that $G\cong G_{n,(g+1)^{r-1}}^{0,\delta}$.

\vspace*{3mm}
\begin{subcase}
$g\leq \delta <g+k$. Note that $g\neq 0$. Otherwise, $\delta< k$, which contradicts $k\leq \delta$.
\end{subcase}
Recall that $u\in V(G)$ such that $d_G(u)=\delta$ and $d_F(u)=t$. First of all, we assert $u\notin F$.
Otherwise, we consider the following two ways according to the value of $k$. \par
\textbf{In the case of $k\geq 2$.} By $u\in F$ and the maximality of $\rho(G)$, using the similar analysis as the proof of Claims \ref{claim5}-\ref{claim7}, we can deduce that $G\cong G_{n,(g+1)^{r-1}}^{k-1,\delta-k+1}$ and $n_1=n-(r-1)(g+1)-k\geq n_2=\cdots= n_r=g+1$.
Let $V(K_{n_i})=\{v_i^1,v_i^2,\ldots,v_i^{n_i}\}$, $P_1=\{v_1^1,v_1^2,\ldots,v_1^{\delta-k+1}\}$ and $N_F(u)=\{w_1,w_2,\ldots, w_{k-1}\}$.
Assume that $\bm{x}$ is the Perron vector of $A(G)$ and $\rho(G)=\rho$.
By symmetry, we may say $x(v_1^i)=x_1$ for $1\leq i\leq \delta-k+1$, $x(v_1^j)=x_1'$ for $\delta-k+2\leq j\leq n_1$ and $x(v)=x_2$ for $v\in V(G)\backslash (F\cup V(K_{n_1}))$. We have $x(w_1)>x_1>x_1'$ by Lemma \ref{lem1}. Moreover,
\begin{equation*}
\begin{aligned}
\rho x_2=(k-1)x(w_1)+gx_2,
\end{aligned}
\end{equation*}
\begin{equation*}
\begin{aligned}
\rho x(u)=(k-1)x(w_1)+(\delta-k+1)x_1,
\end{aligned}
\end{equation*}
\begin{equation*}
\begin{aligned}
\rho x_1'=(k-1)x(w_1)+(\delta-k+1)x_1+(n_1-\delta+k-2)x_1'>(k-1)x(w_1)+(n_1-1)x_1',
\end{aligned}
\end{equation*}
\begin{equation*}
\begin{aligned}
\rho x(w_1)=(k-2)x(w_1)+x(u)+(\delta-k+1)x_1+(n_1-\delta+k-1)x_1'+(r-1)(g+1)x_2,
\end{aligned}
\end{equation*}
from which we have $(\rho-n_1+1)(x_1'-x_2)>(n_1-1-g)x_2\geq 0$. Since $G$ contains $K_{n_1}$ as a proper subgraph, $\rho(G)>\rho(K_{n_1})=n_1-1$. It follows that $x_1'>x_2$. Combining this with $x_1>x_1'$, we get $x_1>x_2$. Then
\begin{align} \label{equ13}
&(n_1-\delta+k-1)(x(u)+x_1')+gx(u)+[(r-2)(g+1)+\delta-k+1]x_2-(k-1-\delta+g)x(w_1)\nonumber\\
=&\frac{(\rho+1)\{(n_1-\delta+k-1)(x(u)+x_1')+gx(u)+[(r-2)(g+1)+\delta-k+1]x_2-(k-1-\delta+g)x(w_1)\}}{\rho+1}\nonumber\\
=&\frac{1}{\rho+1}[(2n_1+(r-2)(g+1))(k-1)x(w_1)+(2n_1+k-\delta-1)(\delta-k+1)x_1+n_1x(u)\nonumber\\
&+(n_1-\delta+k-1)(n_1-g)x_1'+(g+1)(2(r-1)+r(\delta-k)-g)x_2]\nonumber\\
>&\frac{n_1(k-1)x(w_1)+n_1(\delta-k+1)x_1-g(g+1)x_2}{\rho+1}~~(\text{since}~r\geq 2,~n_1\geq g+1~\text{and}~g+k> \delta\geq k)\nonumber\\
>&\frac{n_1\delta x_2-g(g+1)x_2}{\rho+1}~~(\text{since}~x(w_1)>x_1>x_2)\nonumber\\
\geq&0~~(\text{since}~n_1\geq g+1~\text{and}~\delta\geq g).
\end{align}
Note that $\delta-k<g$. Let $E_1=\{uv_1^i \mid \delta-k+2\leq i\leq n_1\}+\{uv_j^l \mid 2\leq j\leq r-1,~1\leq l\leq n_j\}+\{uv_r^h \mid 1\leq h\leq g\}$,
$E_2=\{v_r^{g+1}w_i \mid \delta-g+1\leq i\leq k-1\}$ and $G_1=G-E_2+E_1$. Clearly, $G_1\in\mathcal{G}_n^{k,\delta}$.
Assume that $\bm{y}$ is the Perron vector of $A(G_1)$. By symmetry, we may say that $y(v_1^i)=y_1$ for $1\leq i\leq n_1$, $y(v_j^i)=y_2$ for $2\leq j\leq r-1$ and $1\leq i\leq n_j$, $y(v_r^j)=y_r$ for $1\leq j\leq g$, $y(w_1)=y(w_i)$ for $2\leq i\leq \delta-g$ and $y(w_j)=y(u)$ for $\delta-g+1\leq j\leq k-1$. We have $y(u)>y_1$ and $y_r>y(v_r^{g+1})$ by Lemma \ref{lem1}, then
\begin{equation*}
\begin{aligned}
\rho(G_1)y_1=(\delta-g)y(w_1)+(k-\delta+g)y(u)+(n_1-1)y_1,
\end{aligned}
\end{equation*}
\begin{equation*}
\begin{aligned}
\rho(G_1)y_r=&(\delta-g)y(w_1)+(k-\delta+g)y(u)+(g-1)y_r+y(v_r^{g+1})\\
<&(\delta-g)y(w_1)+(k-\delta+g)y(u)+gy_r,
\end{aligned}
\end{equation*}
from which we have
$$(\rho(G_1)-n_1+1)(y_1-y_r)>(n_1-g-1)y_r\geq 0.$$
Since $G_1$ contains $K_{n_1}$ as a proper subgraph, $\rho(G_1)>\rho(K_{n_1})=n_1-1$. It follows that $y_1>y_r$. Combining this with $y_r>y(v_r^{g+1})$, we have $y_1>y(v_r^{g+1})$. Furthermore,
\begin{equation*}
\begin{aligned}
&\bm{y}^T(\rho(G_1)-\rho(G))\bm{x}\\
=& \bm{y}^T(A(G_1)-A(G))\bm{x}\\
=& \sum_{uv_1^i\in E_1}(x(u)y(v_1^i)+x(v_1^i)y(u))+\sum_{uv_j^l\in E_1}(x(u)y(v_j^l)+x(v_j^l)y(u))+\sum_{uv_r^h\in E_1}(x(u)y(v_r^h)+x(v_r^h)y(u))\\
&-\sum_{v_r^{g+1}w_i \in E_2}(x(v_r^{g+1})y(w_i)+x(w_i)y(v_r^{g+1}))\\
=&(n_1-\delta+k-1)(x(u)y_1+x_1'y(u))+(r-2)(g+1)(x(u)y_2+x_2y(u))+g(x(u)y_r+x_2y(u))\\
&-(k-1-\delta+g)(x_2y(u)+x(w_1)y(v_r^{g+1}))\\
=&(n_1-\delta+k-1)(x(u)y_1+x_1'y(u))+(r-2)(g+1)x(u)y_2+[(r-2)(g+1)+\delta-k+1]x_2y(u)\\
&+gx(u)y_r-(k-1-\delta+g)x(w_1)y(v_r^{g+1})\\
>& y(v_r^{g+1})\left\{(n_1-\delta+k-1)(x(u)+x_1')+gx(u)+[(r-2)(g+1)+\delta-k+1]x_2-(k-1-\delta+g)x(w_1)\right\}\\
~~&(\text{since}~y(u)>y_1>y_r>y(v_r^{g+1})~\text{and}~r\geq 2)\\
>& 0~~(\text{by}~\ref{equ13}),
\end{aligned}
\end{equation*}
and hence $\rho(G_1)>\rho(G)$,
which contradicts the maximality of $\rho(G)$. It is easy to see that $d_{G_1}(v_r^{g+1})=\delta$. This demonstrates that if $k\geq 2$, then $u=v_r^{g+1}\notin F$. \par
\textbf{In the case of $k=1$.} By $u\in F$ and the maximality of $\rho(G)$, using the similar analysis as the proof of Claims \ref{claim6} and \ref{claim7}, we can obtain that $G\cong G_{n,(g+1)^{r-1}}^{0,\delta}$ and $n_1=n-(r-1)(g+1)-1\geq n_2=\cdots=n_r=g+1$. Note that $g\leq \delta<g+1$, then $\delta=g$.
Denote $V(B_i)=\{v_i^1,v_i^2,\ldots,v_i^{n_i}\}$, $P_1=\{v_1^1,v_1^2,\ldots,v_1^{\delta-r+1}\}$ and $P_j=\{v_j^1\}$ for $2\leq j\leq r$. Let $G_2=G+\{uv_1^i \mid \delta-r+2\leq i\leq n_1\}+\{uv_j^l \mid 2\leq j\leq r-1,~1\leq l\leq n_j\}+\{uv_r^h \mid 1\leq h\leq g\}$. Obviously,
$G_2\in\mathcal{G}_n^{k,\delta}$ and $d_{G_2}(v_r^{g+1})=\delta$. According to Lemma \ref{lem2}, we get $\rho(G_2)>\rho(G)$, which contradicts the maximality of $\rho(G)$. This means that if $k=1$, then $u=v_r^{g+1}\notin F$. \par
Based on the above results, we can conclude that if $g\leq \delta <g+k$, then $u \notin F$. Next,
using the similar analysis as the proof of Claims \ref{claim2}-\ref{claim4}, we can deduce that $G\cong G_{n,(g+1)^{r-1}}^{\delta-g,g}$, as desire.

\vspace*{3mm}
\begin{subcase}
$\delta\geq g+k$.
\end{subcase}
Recall that $|V(B_i)|=n_i$ for $1\leq i\leq q$ and $q=r$. Therefore,
$$\rho(G)\leq \rho(K_k\vee(K_{n_1}\cup K_{n_2}\cup \cdots \cup K_{n_r})),$$
with equality if and only if $G\cong K_k\vee(K_{n_1}\cup K_{n_2}\cup \cdots \cup K_{n_r})$. Without loss of generality,
we may assume that $n_1\geq n_2\geq \cdots\geq n_r$. Obviously, $n_r\geq \delta-k+1$ because the minimum degree of $G$ is $\delta$. Combining this with Lemma \ref{lem4}, we get
$$\rho(K_k\vee(K_{n_1}\cup K_{n_2}\cup \cdots \cup K_{n_r}))\leq \rho(K_k\vee(K_{n-k-(\delta-k+1)(r-1)}\cup (r-1)K_{\delta-k+1})),$$
with equality if and only if $(n_1,n_2,\ldots,n_r)=(n-k-(\delta-k+1)(r-1),\delta-k+1,\ldots,\delta-k+1)$. By the maximality of $\rho(G)$, we conclude that $G\cong K_k\vee(K_{n-k-(\delta-k+1)(r-1)}\cup (r-1)K_{\delta-k+1})$.\par
This completes the proof.
$\hfill\qedsymbol$

\section{The maximum spectral radius of graphs with $g$-good neighbor connectivity} \label{sec3}
Recall that $F\subseteq V(G)$ is called a $g$-good neighbor cut of $G$ if $G-F$ is disconnected and each vertex $v\in V(G)\backslash F$ has at least $g$ neighbors. The $g$-good neighbor connectivity of $G$, denoted by $\kappa_g(G)$ or $\kappa_g$, is the minimum cardinality of $g$-good neighbor cuts of $G$. Obviously, $\kappa_0(G)=\kappa(G)$.
As a conditional connectivity, the $g$-good neighbor connectivity has improved networks' fault tolerance and has been extensively employed in the analysis of various networks, see \cite{H.Q. Liu, Q.R. Zhou}.
Similar to $g$-good $r$-component connectivity, based on the Brualdi-Solheid problem, it is natural for us to
consider the following problem.
\begin{problem}\label{problem3}
Which graphs attain the maximum spectral radii among all connected graphs of order $n$ with fixed minimum degree $\delta$ and $g$-good neighbor connectivity $\kappa_g(G)$?
\end{problem}

Let $\mathcal{G}_n^{\kappa_g,\delta}$ be the set of graphs of order $n$ with minimum degree $\delta$ and $g$-good neighbor connectivity $\kappa_g$. The condition of $g$-good neighbor connectivity is weaker than that of $g$-good $r$-component connectivity, and there is no requirement on the number of components.
Therefore, we know that the graphs which achieve the maximum spectral radii among all graphs in $\mathcal{G}_n^{\kappa_g,\delta}$ must be in $\mathcal{G}_n^{k,\delta}$. Furthermore, we have $r=2$ by Lemma \ref{lem2}.
Let $G_{n,g+1}^{\delta,0}$ be the graph obtained from $K_1\cup(K_{\kappa_g-1}\vee(K_{n-g-1-\kappa_g}\cup K_{g+1}))$ by adding $\delta$ edges between the isolated vertex $K_1$ and $K_{\kappa_g-1}$.
Let $G_{n,g+1}^{\delta-g,g}$ be the graph obtained from $K_1\cup(K_{\kappa_g}\vee(K_{n-g-1-\kappa_g}\cup K_g))$ by adding $\delta-g$ edges between the isolated vertex $K_1$ and $K_{\kappa_g}$, and adding $g$ edges between $K_1$ and $K_g$.
Let $G_{n,g+1}^{\kappa_g-1,\delta-\kappa_g+1}$ be the graph obtained from $K_1\cup (K_{\kappa_g-1}\vee(K_{n-g-1-\kappa_g}\cup K_{g+1}))$ by adding $\kappa_g-1$ edges between the isolated vertex $K_1$ and $K_{\kappa_g-1}$, and adding $\delta-\kappa_g+1$ edges between $K_1$ and $K_{n-g-1-\kappa_g}$.
Let $G_{n,g+1}^{0,\delta}$ be the graph obtained from $K_1\cup (K_{n-g-2}\cup K_{g+1})$ by adding $\delta-1$ edges between the isolated vertex $K_1$ and $K_{n-g-2}$, and adding one edge between $K_1$ and $K_{g+1}$, respectively.
Then we have the following result.

\begin{thm}
Let $G\in \mathcal{G}_n^{\kappa_g,\delta}$, where $n\geq \kappa_g+2(g+1)$. Then we have the following statements.
 \begin{description}
   \item[(I)] If $\kappa_g>\delta$ and $\delta< g$, then $\rho(G)\leq \rho(G_{n,g+1}^{\delta,0})$ with equality if and only if $G\cong G_{n,g+1}^{\delta,0}$;
   \item[(II)] If $\kappa_g>\delta\geq g$, then $\rho(G)\leq \rho(G_{n,g+1}^{\delta-g,g})$ with equality if and only if $G\cong G_{n,g+1}^{\delta-g,g}$;
   \item[(III)] If $2\leq\kappa_g\leq\delta<g$, then $\rho(G)\leq \rho(G_{n,g+1}^{\kappa_g-1,\delta-\kappa_g+1})$ with equality if and only if $G\cong G_{n,g+1}^{\kappa_g-1,\delta-\kappa_g+1}$;
   \item[(IV)] If $1=\kappa_g\leq\delta<g$, then $\rho(G)\leq \rho(G_{n,g+1}^{0,\delta})$ with equality if and only if $G\cong G_{n,g+1}^{0,\delta}$;
   \item[(V)] If $\kappa_g\leq\delta$ and $g\leq \delta <g+\kappa_g$, then $\rho(G)\leq \rho(G_{n,g+1}^{\delta-g,g})$ with equality if and only if $G\cong G_{n,g+1}^{\delta-g,g}$;
   \item[(VI)] If $\delta \geq g+\kappa_g$, then $\rho(G)\leq \rho(K_{\kappa_g}\vee(K_{n-\delta-1}\cup K_{\delta-\kappa_g+1}))$ with equality if and only if $G\cong K_{\kappa_g}\vee(K_{n-\delta-1}\cup K_{\delta-\kappa_g+1})$.
 \end{description}
\end{thm}

\end{document}